\documentclass[12pt]{amsart}
\input{amssymb.sty}

\newtheorem{theorem}{Theorem}

\newtheorem{lemma}[theorem]{Lemma}
\newtheorem{definition}[theorem]{Definition}
\newtheorem{corollary}[theorem]{Corollary}
\newtheorem{proposition}[theorem]{Proposition}

\newtheorem{example}[theorem]{Example}

\begin{document}
\title[Residual Solvability]
{On the residual solvability of generalized free products of
finitely generated nilpotent groups}
\author[D.Kahrobaei]{D.Kahrobaei}
\address{Mathematical Institute, University of St
Andrews, North Haugh, St Andrews, Fife, KY16 9SS, UK}
\email{delaram.kahrobaei@st-andrews.ac.uk}
\urladdr{http://www-groups.mcs.st-and.ac.uk/${\sim}$delaram/}

\begin{abstract} In this paper we study the residual
solvability of the generalized free product of finitely generated
nilpotent groups. We show that these kinds of structures are often
residually solvable.
\end{abstract}
\thanks{The research of this author has been supported by CUNY research
foundation at the City College of the City University of New York
and New York Group Theory Cooperative.}

\date{\today} 


\dedicatory{Dedicated to Gilbert Baumslag on the occasion of his
70th Birthday}

\maketitle
\section{Introduction and Motivation}
The notion of residual properties was first introduced by Philip
Hall in 1954 \cite{PH54}. Let $X$ be a class of groups. $G$ is
residually-$X$ if, for every non-identity element $g$ in $G$,
there is an epimorph of $G$ in $X$ such that the element
corresponding to $g$ is not the identity.

In this paper we focus on residual solvability. The notion of
residual solvability of groups was pioneered particularly by
Gilbert Baumslag in his celebrated paper \cite{GB71} where he
proved that positive one-relator groups are residually solvable.
The author in \cite{DK} studied this notion in general. In
\cite{DK2005} the author studied residual solvability of
residually solvable groups.

In 1963 Gilbert Baumslag studied residual finiteness of
generalized free products of finitely generated torsion-free
nilpotent groups \cite{GB63}. He showed first that if both factors
are non-abelian then this amalgamation is not residually-finite
and he found conditions that ensure $G$ is residually finite. He
actually proved in general that this kind of structure is
free-by-residually-finite, or meta-residually-finite. A few years
later, in 1968, one of his students, Joan Landman Dyer, continued
his work \cite{JD68} on residual finiteness of these structures.
She showed in particular that if the factors are not finitely
generated nilpotent groups, these kinds of groups are not
free-by-residually-finite (by taking isomorphic factors of class
three), but rather are
residually-finite-by-free-by-residually-finite.

It is interesting to mention that the generalized free product of
two finitely generated nilpotent groups (or two finitely generated
free groups) with a finitely generated subgroup amalgamated has a
solvable word problem (see page 150 \cite{GB93}). In this work we
consider the generalized free product of finitely generated
nilpotent groups, and discuss how close such groups are to being
residual solvable. We show how conditions on the amalgamating
subgroup affects residual solvability. Note that the
free products of finitely generated nilpotent groups are residually solvable.\\
\newline
\noindent {\bf Effect of abelianization on amalgamated products of
nilpotent groups}
\newline
\noindent We give a complete explanation of abelianization of
generalized free products. In particular, we see how the
abelianization of such groups with nilpotent factors is not
trivial, using the Frattini theorem, which states that the
commutator subgroup of a finitely generated nilpotent group is a
non-generator set. As a result, we conclude that generalized free
products of two finitely generated nilpotent groups are not
perfect. Note that the condition of being finitely generated is
required in these results since we use the Frattini theorem.
\begin{theorem} \label{Not_Perfect} The generalized free product
of two finitely generated nilpotent groups amalgamated by a proper
subgroup of them is not perfect.
\end{theorem}
\noindent This property of finitely generated nilpotent groups
plays a critical role in our approach to our question. For example
finitely generated polycyclic factors that satisfy the maximal
condition, but fail to satisfy the condition of the Frattini
theorem, turn out to be perfect under certain conditions.\\\\
\noindent
{\bf Cyclic amalgamated subgroup and residual solvability}
\newline
\noindent We first consider the case when the amalgamating
subgroup is cyclic. It turns out that choosing an appropriate
solvable filtration of each factor so that the generator of
amalgamating subgroup does not lie in the $n$-th term of upper
central series enables us to prove that it is residually
solvable-by-solvable.
\begin{theorem} \label{Cyclic_amalgam} The generalized free product of two
finitely generated nilpotent groups amalgamated by a cyclic
subgroup is residually solvable.
\end{theorem}
\noindent {\bf The amalgamated subgroup is central in all factors}
\newline
\noindent We next define the generalized central product of an
arbitrary number of groups, and show that each factor injects into
such a product. We then look at the case where the amalgamating
subgroup is central in all factors, and we show using the
generalized central products that these structures turn out to be
free-by-nilpotent. We note that a nilpotent extension of a free
group is not residually nilpotent but rather is residually
solvable.
\begin{theorem} \label{finiteNilpotentRS} The generalized free product of an arbitrary number of finitely
generated nilpotent groups of bounded class, amalgamating a
central subgroup in each of the factors, is an extension of a free
group by a nilpotent group. Furthermore, such groups are
residually solvable.
\end{theorem}
\noindent {\bf The case where one of the factors is abelian}
\newline
\noindent We then consider the case where just one of the factors
is abelian. It turns out that such groups with an abelian factor
are free-by-nilpotent.
\begin{theorem} \label{abelianonefactor} The generalized free product of a finitely generated
torsion-free abelian and a nilpotent group is residually
solvable-by-abelian-by-finite abelian. Furthermore such groups are
residually solvable.
\end{theorem}
\noindent {\bf The case when the amalgamated subgroup is of finite
index in at least one of the factors}
\begin{theorem} \label{finind} The generalized free product of two finitely
generated torsion-free nilpotent groups amalgamating a proper
subgroup which is of finite index in at least one of the factors,
is an extension of a free group by a torsion-free nilpotent group.
Furthermore, such groups are residually solvable.
\end{theorem}
Note that here the condition of being torsion-free is necessary
for the factors, because of Mal'cev's Fundamental Theorem, as we
see later.\\\\
\noindent {\bf Doubles of nilpotent groups and residual
solvability, arbitrary number of factors}
\newline
\noindent By a double we mean the amalgamated product of two
groups where the factors are isomorphic and the amalgamated
subgroups are identified under the same isomorphism. One can
generalize the definition of doubles to arbitrary number of
factors.
\begin{theorem}\label{finite_double}
Let $\{A_i| i \in I \}$ be an arbitrary indexed family of
isomorphic torsion-free nilpotent groups, such that $\bigcap_{i
\in I} A_i = C$, and let $G$ be the generalized free product of
$A_i$s amalgamated by $C$:
\begin{eqnarray*}
G = \{ {\prod_{i \in I}}^* A_i; C \}.
\end{eqnarray*}
Then $G$ is an extension of a free by nilpotent group.
Furthermore, $G$ is residually solvable.
\end{theorem}
\noindent {\bf Generalized free products of nilpotent groups
sometimes may fail to be residually solvable.}
\newline
\noindent Gilbert Baumslag gave a counter example in \cite{GB72}
that shows that not every subgroup of the generalized free product
of two finitely generated torsion-free nilpotent groups is
indicable. We use this example to prove the following proposition:
\begin{proposition} \label{nil_neg} There exist two finitely generated torsion-free non-abelian free
nilpotent groups $A$ and $B$ such that the amalgamated product of
them with abelian amalgamating subgroup $C_A = C_B$,
\begin{eqnarray*}
G = \{ A \ast B ; C_A = C_B \},
\end{eqnarray*}
is not residually solvable and not poly-residually solvable.
\end{proposition}
This implies that an abelian amalgamating subgroup is not
sufficient for residual solvability.\\
\newline
\noindent {\bf Poly-residual solvability of the generalized free
product of nilpotent groups}\\
\noindent As Proposition \ref{nil_neg} suggests, we need to impose
some conditions on the amalgamating subgroup to ensure that the
generalized free product of finitely generated nilpotent group be
poly-residually solvable; by a poly-residually solvable group we
mean a group that has at least one poly-residually solvable series
(see Section \ref{Poly_RS} for precise definition). Here is the
theorem that gives us these conditions:
\begin{theorem} \label{Poly-RS} The generalized free product of two finitely generated
nilpotent groups, $A$ and $B$ amalgamating subgroups of them,
$C_A$, $C_B$ respectively is poly-residually solvable if the
following condition holds: the solvable filtrations of $A$ and $B$
which are defined on the upper central series of $A$ and $B$
respectively, are solvable with compatible filtrations, i.e.
$\xi_i A \cap C_A \stackrel{\phi_i}{=} \xi_i B \cap C_B$.
\end{theorem}

\subsection*{Acknowledgment} I thank my Ph.D. supervisor
G.Baumslag and also K.J.Falconer and P.de la Harpe for helpful
comments.
\section{Background and Preliminary Results}
In this section we recall some definitions and facts and prove
some lemmas to be used later.

Recall $G$ is called an extension of $A$ by $Q$ if there exist $A$
and $Q$ and a short exact sequence $1 \rightarrow A \rightarrow G
\rightarrow Q \rightarrow 1$. We say $G$ is meta-$X$, where $X$ is
a property (or class), if $G$ is an extension of $A$ by $Q$ where
$A$ and $Q$ have property (or class) $X$.
\subsection{Subgroups of amalgamated products}
\label{Hanna} We will use a theorem of Hanna Neumann \cite{HN49}
extensively in this paper. With regard to abstract groups, Hanna
Neumann showed in the 1950s that, in general, subgroups of
amalgamated products are no longer amalgamated products, but
generalized free products, indeed she proved the following: let
$K$ be a subgroup of $G =\{ A \ast B; C \}$, then $K$ is an
HNN-extension of a tree product in which the vertex groups are
conjugates of subgroups of either $A$ or $B$ and the edge groups
are conjugates of subgroups of $C$. The associated subgroups
involved in the HNN-extension are also conjugates of subgroups of
$C$. If $K$ misses the factors $A$ and $B$ (i.e. $K \cap A = \{1\}
= K \cap B$), then $K$ is free; and if $K$ misses the amalgamated
subgroup $C$ (i.e. $K \cap C = \{1\}$), then $K = {\prod_{i \in
I}}^* X_i \ast F$, where the $X_i$ are conjugates of subgroups of
$A$ and $B$ and $F$ is free (see \cite{GB93} for more
information).

Let us mention that later a description was given by the
Bass-Serre theory \cite{S80}, with groups acting on graphs to give
geometric intuition: the fundamental group of a graph of groups
generalizes both amalgamated products, HNN-extensions and tree
products.
\section{Some results on the structure of the abelianization}
\label{Abelian} The following lemma formulates the abelianization
of the amalgamated products of two groups. (We leave it as an
exercise to the reader to check this.)
\begin{lemma}
\label{abelianization_lemma} Let $G$ be the amalgamated product of
two groups $A$ and $B$ amalgamating $C$,
\begin{eqnarray*}
G=\{A \ast B; C \}.
\end{eqnarray*}
Then the abelianization of $G_{ab}$ of $G$ takes the following
form:
\begin{eqnarray*}
G_{ab} = (A_{ab} \times B_{ab})/{gp(\bar{c} \alpha \bar{c}^{-1}
\beta | c \in C)},
\end{eqnarray*}
where $\alpha$ is the monomorphism from $C$ into $A$ and $\beta$
is the monomorphism from $C$ into $B$.
\end{lemma}
\begin{lemma} \label{onto_abelianization} Let $G$ be the amalgamated product of $A$ and $B$ identifying $C_A$ with $C_B$,
\begin{eqnarray*}
G = \{A \ast B;C_A = C_B\}.
\end{eqnarray*}
Then the abelianization $G_{ab}$ of $G$ maps onto
\begin{eqnarray*}
D = A_{ab}/{{gp({\bar{c}}_a | c_a \in C_A)}} \times
B_{ab}/{gp({{\bar{c}}_b}^{-1} | c_b \in C_B)}.
\end{eqnarray*}
\end{lemma}
\begin{proof}
By Lemma \ref{abelianization_lemma}, the abelianization of $G$ can
be expressed as:
\begin{eqnarray*}
G_{ab} = (A_{ab} \times B_{ab})/{gp(\bar{c_a}{\bar{c_b}}^{-1})}.
\end{eqnarray*}
Put
\begin{eqnarray*}
       N   & = & gp(\bar{c_a}{\bar{c_b}}^{-1})\\
       N_1 & = & {gp({\bar{c}}_a)}\\
       N_2 & = & {gp({{\bar{c}}_b}^{-1})}.
\end{eqnarray*}
Let $\theta$ be the map from $G_{ab}$ into $D$,
\begin{eqnarray*}
     \theta :  G_{ab} & \rightarrow & D\\
     \theta : (A_{ab} \times B_{ab})/N & \rightarrow & A_{ab}/{N_1} \times B_{ab}/{N_2}
\end{eqnarray*}
defined by:
\begin{eqnarray*}
\bar{a} N \mapsto \bar{a} N_1 \text{ and } \bar{b} N \mapsto
\bar{b} N_2.
\end{eqnarray*}
A typical element $w$ in $N$ takes the form $w = \bar{c_a}
{\bar{c_b}}^{-1}$, so
\begin{eqnarray*}
w = \bar{c_a} {\bar{c_b}}^{-1} \mapsto (\bar{c_a} N_1,
{\bar{c_b}}^{-1} N_2) = (N_1 , N_2).
\end{eqnarray*}
and by Von Dyck's theorem $\theta$ is an epimorphism as requested.
\end{proof}
The following is an alternative proof of Lemma
\ref{onto_abelianization}, using presentations.
\begin{proof}
Let $A$ and $B$ have presentations:
\begin{eqnarray*}
A = \langle X;R \rangle  \; \text{ and } B = \langle Y;S \rangle.
\end{eqnarray*}
So the abelianizations of $A$ and $B$ have presentations:
\begin{eqnarray*}
              A_{ab} & = & \langle X; R \cup \{[x_i, x_j]| x_i, x_j \in X\} \rangle \\
          B_{ab} & = & \langle Y; S \cup \{[y_i, y_j]| y_i, y_j \in Y\} \rangle .
\end{eqnarray*}
The direct product of $A_{ab}$ and $B_{ab}$ has the following
presentation since, by definition of the direct product, all
elements of $A_{ab}$ commute with all elements of $B_{ab}$:
\begin{eqnarray*}
A_{ab} \times B_{ab} = \langle X \cup Y;&& R \cup S \cup \{  [x_i, x_j]|\; x_i, x_j \in X \}\\
                                 && \cup \{  [y_i, y_j]|\; y_i, y_j \in Y \}\\
                 && \cup \{\;[x, y]  \;|\; x \in X, y \in Y \} \rangle .
\end{eqnarray*}
With $\alpha$ the monomorphism from $C$ into $A$, and $\beta$ the
monomorphism from $C$ into $B$, put
\begin{eqnarray*}
N   & = & gp(c \alpha \; c {\beta}^{-1} | c \in C )\\
N_1 & = & gp(c \alpha | c \in C)\\
N_2 & = & gp(c {\beta}^{-1} | c \in C).
\end{eqnarray*}
By Lemma \ref{abelianization_lemma}, $G_{ab}$, the abelianization
of $G$, has the following presentation:
\begin{eqnarray*}
G_{ab} = (A_{ab} \times B_{ab})/N =  \langle X \cup Y; R \cup S & \cup & \{[x_i, x_j]| x_i, x_j \in X\} \\
                                                         & \cup & \{[y_i, y_j]| y_i, y_j \in Y\} \\
                             & \cup & \{[x, y]\;\;| x \in X, y \in Y \} \\
                             & \cup & \{ c \alpha \; c {\beta}^{-1} | c \in C \} \rangle .
\end{eqnarray*}
Then the presentation of each factor of $D$ is
\begin{eqnarray*}
{A_{ab}}/{N_1}  & = & \langle X; R \cup \{ [x_i, x_j]|\; x_i, x_j \in X \} \cup \{c \alpha       | c \in C\} \rangle \\
{B_{ab}}/{N_2}  & = & \langle Y; S \cup \{ [y_i, y_j]|\; y_i, y_j
\in Y \} \cup \{c {\beta}^{-1} | c \in C\} \rangle .
\end{eqnarray*}
Their direct product $D$ has the following presentation:
\begin{eqnarray*}
D = {A_{ab}}/{N_1} \times {B_{ab}}/{N_2} = && \langle  X \cup Y; R \cup S \cup \{ [x_i, x_j]|\; x_i, x_j \in X \}\\
                                           && \cup \{ [y_i, y_j]|\; y_i, y_j \in Y \}\\
                       && \cup \{ c \alpha| c \in C\}\\
                       && \cup \{ c {\beta}^{-1} | c \in C\}\\
                       && \cup \{ [x ( c \alpha), y (c {\beta}^{-1})]| x \in X, y \in Y, c \in C \} \rangle .
\end{eqnarray*}
Now define $\theta$ to be a map from $G_{ab}$ into $D$, by sending
$x_i \mapsto x_i$ and $y_i \mapsto y_i$. Then $N \mapsto 1$, so by
Von Dyck's Theorem, $\theta$ defines a homomorphism from $G_{ab}$
onto $D$.
\end{proof}
\section{Effects of the order of the abelianization of amalgamated products of nilpotent groups}
In order to study the effect of the order of the abelianization of
the amalgamated products of nilpotent groups, we first study the
effect of indicability on the order of finitely generated
nilpotent groups. Recall that a group $A$ is termed indicable if
there exists a homomorphism of $A$ onto the infinite cyclic group.
A finitely generated group $G$ is indicable if and only if
$G_{ab}$ is infinite. Higman \cite{GH40} proves that every
finitely generated torsion-free nilpotent group is indicable and
is therefore infinite. We conclude that the abelianization of
every finitely generated torsion-free nilpotent group is also
infinite. Note that the abelianization of an infinite finitely
generated nilpotent group $A$ is again infinite, since there is a
canonical homomorphism from $A$ onto $A/{\tau A}$, where $\tau A$
is the torsion group of $A$. Since $A/{\tau A}$ is a finitely
generated torsion-free nilpotent group, $A/{\tau A}$ is infinite.
Hence $A$ is infinite.

To prove the following proposition we use a theorem of G. Baumslag
\cite{GB72} that states that the amalgamated product of two
finitely generated torsion-free nilpotent groups is indicable.
\begin{proposition} \label{Abel_Amal_t-f_Nilp} The abelianization of the amalgamated product of two
finitely generated nilpotent groups is infinite only if one of the factors
is infinite.
\end{proposition}
\begin{proof} First let $A$ and $B$ be two finitely generated torsion-free nilpotent
groups, let $C$ be a proper subgroup of both of them, and $G$ be
the amalgamated product of them, amalgamating $C$,
\begin{eqnarray*}
G= \{A \ast B;C\}.
\end{eqnarray*}
$G$ is indicable by the theorem of G.Baumslag. Therefore $G_{ab}$,
the abelianization of $G$, is infinite (see \cite{GB93} for more
information). Now let just one of the factors be torsion-free.
There exists an epimorphism from $G$ onto a finitely generated
torsion-free nilpotent group,
\begin{eqnarray*}
G^* = \{A/{\tau A} \ast B/{\tau B} ; C/{\tau C}\},
\end{eqnarray*}
the amalgamated product of finitely generated torsion-free
nilpotent groups. Therefore ${G^*}_{ab}$ is infinite, using the
first case, Implying that $G_{ab}$ is infinite.
\end{proof}
\section{Proof of Theorem \ref{Not_Perfect}}
In this section we first recall the Frattini Subgroup Theorem.
Recall that the Frattini subgroup of $G$, $\Phi G$, is the
intersection of all maximal subgroups of $G$. In particular, $\Phi
G$ is characteristic. We also recall that $g \in G$ is called a
non-generator of $G$, if whenever $G=gp(X,g)$ then $G=gp(X)$. A
theorem of Frattini states that $\Phi G$ is the set of all
non-generators in $G$. In particular if $G$ is a finitely
generated nilpotent group, then $\delta_2 G \leqslant \Phi G$,
where $\delta_2 G = [G, G]$ (see \cite{DR95} for more
information).
\begin{proof} Let $A$ and $B$ be two finitely generated nilpotent
groups, and let $C_A$ and $C_B$ be proper subgroups of $A$ and $B$
respectively. We want to show that the generalized free product $G
= \{ A \ast B ; C_A = C_B \}$ is not perfect, that is that the
abelianization $G_{ab}$ of $G$ is not trivial.
\newline
\noindent By Lemma \ref{onto_abelianization}, $G_{ab}$ maps onto
\begin{eqnarray*}
D = A_{ab}/{{gp({\bar{c}}_a | c_a \in C_A)}} \times
B_{ab}/{gp({{\bar{c}}_b}^{-1} | c_b \in C_B)}.
\end{eqnarray*}
\noindent
where $\bar{c_a}$ and $\bar{c_b}$ are the images of $c_a
\in C_A$ and $c_b \in C_B$ in $A_{ab}$ and $B_{ab}$ respectively.
\newline
We claim that under the conditions of the Theorem $D$ is not
trivial. In order to prove the claim, we will show that
\begin{eqnarray*}
A_{ab}/{gp({\bar{c}}_a)} \not= \{1\}  \;\; \text{ and }
B_{ab}/{gp({{\bar{c}}_b}^{-1})} \not= \{1\}.
\end{eqnarray*}
Since $C_A$ is a proper subgroup of $A$, the Frattini Subgroup
Theorem implies that $gp(C_A , \delta_2 A)$ is a proper subgroup
of $A$, where $\delta_2 A= [A,A]$. Since $C_A$ is a proper
subgroup of $A$, $gp(C_A , \delta_2 A)$ is also a proper subgroup
of $A$, so
\begin{eqnarray*}
A/{gp(C_A , \delta_2 A)} & \not\simeq & \{1\}.
\end{eqnarray*}
\noindent
By using the third isomorphism theorem we have:
\begin{eqnarray*}
\frac{A_{ab}}{gp(\bar{c_a})} \simeq \frac{\frac{A}{\delta_2
A}}{\frac{gp(C_A, \delta_2 A)}{\delta_2 A}} \simeq \frac{A}{gp(C_A
,\delta_2 A)} \not\simeq \{1\}.
\end{eqnarray*}
\noindent Similarly
\begin{eqnarray*}
{B_{ab}} /{gp({\bar{c_b}}^{-1})} \not= \{1\}.
\end{eqnarray*}
\noindent This proves our claim that, $D \not= \{1\}$ and so
$G_{ab} \not= \{1\}$. Hence $G$ is not perfect.
\end{proof}
\section{Proof of Theorem \ref{Cyclic_amalgam}}
\begin{proof} Let $A$ and $B$ be two finitely generated nilpotent
groups. Let $a$ be a non-identity element of $A$, and $b$ be a
non-identity element of $B$. We can find $m \geqslant 1$ and $n
\geqslant 1$ such that $1 \not= a \in {\xi}_{m+1} A \backslash
{\xi}_{m} A$, and $1 \not= b \in {\xi}_{n+1} B \backslash
{\xi}_{n} B$, (where ${\xi}_i A$ is the $i-$th term of the upper
central series of $A$, and ${\xi}_j B$ is the $j-$th term of the
upper central series of $B$). Let $G=\{A \ast B; a=b \}$ be the
generalized free product of $A$ and $B$ amalgamating $a$ with $b$.
Let $D$ be the central product of $A/{{\xi}_m A}$ and $B/{{\xi}_n
B}$ amalgamating $a{\xi}_m A$ with $b{\xi}_n B$,
\begin{eqnarray*}
D= \{{A/{{\xi}_m A}} \times B/{{\xi}_n B}; a{\xi}_m A= b{\xi}_n B
\}.
\end{eqnarray*}
Note that
\begin{eqnarray*}
C {{\xi}_{n}}B / {{\xi}_{n}}B \simeq C / {C \cap {{\xi}_{n}}}B
\end{eqnarray*}
and
\begin{eqnarray*}
C {{\xi}_{m}}A / {{\xi}_{m}}A \simeq C / {C \cap {{\xi}_{m}}}A
\end{eqnarray*}
are cyclic groups. This confirms the fact that
\begin{eqnarray*}
a{\xi}_m A= b{\xi}_n B.
\end{eqnarray*}
Map $G$ into $D$ and let $K$ be the kernel of this map. Then
observe that $K \cap C = \{1\}$, where $C = gp(a)=gp(b)$.
Therefore $K$ is a free product of conjugates of subgroups of $A$
and $B$, and a free group by the Hanna Neumann Theorem, see
Section \ref{Hanna}. So $K$ is residually solvable. Hence $G$ is
an extension of a residually solvable group by a solvable group.
\end{proof}
\section{The amalgamated subgroup is central in all factors}
\subsection{Generalized central products and some related results}
Here we introduce a general definition of the generalized central
product of an arbitrary number of factors. Let us mention the work
of D.Robinson in \cite{DR95} for finite number of factors. We
adopt the notation used for generalized free products \cite{GB93}:
\begin{definition} \label{gen_cent_prod_def}
\end{definition}
\noindent
Suppose that
\begin{eqnarray*}
\{ A_i = \langle X_i;R_i \rangle  |i \in I\}
\end{eqnarray*}
is an indexed family of presentations
\begin{eqnarray*}
A_i = \langle  X_i ; R_i  \rangle
\end{eqnarray*}
of the groups $A_i$, and suppose $C$ is another group equipped
with monomorphisms
\begin{eqnarray*}
\phi_i : C \rightarrow A_i \; \text{ and } C \leqslant \xi A_i
\;\; (\text{for all } i \in I, \text{ where } \xi A_i \text{ is
the center of } A_i).
\end{eqnarray*}
We term the group $A$ defined by the presentation
\begin{eqnarray*}
A = \langle  \cup X_i ; \cup R_i \cup \{ c \phi_i c^{-1} \phi_j |
c\in C, i, j \in I \} \cup \{ [x_i,x_j]=1 \; i,j \in I \} \rangle,
\end{eqnarray*}
where we assume that the $X_i$ are disjoint, the generalized
central product of $A_i$ amalgamating the central subgroup $C$.
\begin{eqnarray*}
A = {\prod_{i \in I}}^{\times} \{A_i ; C\}.
\end{eqnarray*}
If $C=1$, then $A$ is termed the direct product of the $A_i$.\\
According to Von Dyck's theorem, there are canonical homomorphisms
$\mu_i$ of each $A_i$ to $A$. We will prove that the $\mu_i$ are
monomorphisms, and, if we identify $A_i$ with $A_i \mu_i$, then
\begin{eqnarray*}
c \phi_i = c \phi_j \; \text{for all } c \in C, i,j \in I.
\end{eqnarray*}
Thus, we can identify $C$ with any of its images $C \phi_i$ which
are already identified with $C \phi_j \mu_i$. Then $A_i \cap A_j =
C \; (i,j \in I, i \not= j)$ and $A=gp({\cup}_{i \in I}A_i)$. $A$
can also be written as
\begin{eqnarray*}
A = {({\prod_{i \in I}}^{\times} A_i)}/{gp(c \phi_i c^{-1} \phi_j|
i,j \in I,\; c \in C)}.
\end{eqnarray*}
\begin{lemma} \label{central_any_factor}
$\mu_i : A_i \rightarrow {({\prod_{i \in I}}^{\times} A_i)} / gp(c
\phi_i c^{-1} \phi_j | i, j \in I,  c\in C)$ is a monomorphism for
all $i \in I$.
\end{lemma}
\begin{proof}
Put $S = gp(c \phi_i c^{-1} \phi_j | i, j \in I,  c \in C),
\;\text{ where } C \stackrel{\phi_i}{\rightarrow} A_i.$ To show
that $\mu_i$ is a monomorphism we must show that $\ker \mu_i =1$.
Let $a \in A_i$ and $a \in \ker \mu_i$, so that $a \mu_i = 1$. We
want to show that $a=1$. A generic element in $S$ has the
following form:
\begin{eqnarray*}
s = (c_1 \phi_{i_1} {c_1}^{-1} \phi_{j_1}) \cdots (c_n \phi_{i_n}
{c_n}^{-1} \phi_{j_n}).
\end{eqnarray*}
Since $a \mu_i = 1$, then $a \in S$, and
\begin{eqnarray*}
a = (c_1 \phi_{i_1} {c_1}^{-1} \phi_{j_1}) \cdots (c_n \phi_{i_n}
{c_n}^{-1} \phi_{j_n}).
\end{eqnarray*}
Let us consider two cases: the case where none of $i_r$ and $j_r$
are equal to $i$,
and the case where some of them are equal to $i$.\\
Case 1: Since none of the subscripts are equal to $i$, this
implies that
\begin{eqnarray*}
a \in gp(A_k |k \not=i).
\end{eqnarray*}
But $a$ is also an element in $A_i$. Therefore
\begin{eqnarray*}
a \in A_i \cap gp(A_k | i \not= k) =1
\end{eqnarray*}
(by a property of the direct product), so $a=1$.
\newline
Case 2: Now suppose that some of the indexes are equal to $i$, say
$i_l =i$. Note that for $c \in C$, $c \phi_i$ is central in $A_i$.
Thus we have:
\begin{eqnarray*}
a=(\prod_{i_l} c_l \phi_{i_l} {c_l}^{-1} \phi_{j_l})(\prod_{k_m
\not= i} c_1 \phi_{k_1} {c_1}^{-1} \phi_{k_1} \cdots c_n
\phi_{k_n} {c_n}^{-1} \phi_{k_n}).
\end{eqnarray*}
By a similar argument,
\begin{eqnarray*} (\prod_{k_m \not= i} c_1 \phi_{k_1}
{c_1}^{-1} \phi_{k_1} \cdots c_n \phi_{k_n} {c_n}^{-1} \phi_{k_n}
) = 1.
\end{eqnarray*}
Now
\begin{eqnarray*}
\prod_{i_l} c_l \phi_{i_l} & = & ({\prod} c_l) \phi_{i} = c
\phi_i.
\end{eqnarray*}
Hence
\begin{eqnarray*}
                         a & = & ( c \phi_{i})(\prod_{j_l \not= i} {c_l}^{-1} \phi_{j_l}).
\end{eqnarray*}
We have
\begin{eqnarray*}
a c^{-1} \phi_{i} = (\prod_{j_l \not= i} {c_l}^{-1} \phi_{j_l})
\in A_i \cap gp(A_k| k \not=i) =1,
\end{eqnarray*}
so $a = c \phi_{i}$. This implies that $a= 1$, since
\begin{eqnarray*}
a= c \phi_i = (\prod_{i_l} c_l \phi_{i_l} {c_l}^{-1} \phi_{j_l}).
\end{eqnarray*}
The right hand side has an even number of factors, so this is
possible only if it is equal to $1$, as required.
\end{proof}
\subsection{Proof of Theorem \ref{finiteNilpotentRS}}
Before proving the Theorem, we bring together some of the facts
and related lemmas that we will use in the proof. One of these
facts is that a direct product of finitely many nilpotent groups
is nilpotent. However we note that a direct product of infinitely
many nilpotent groups is not necessarily nilpotent. (see D.Segal's
book \cite{DS83} for more information)
\begin{lemma} \label{Central_Prod_Nil} The generalized central product of finitely many nilpotent groups is nilpotent.
\end{lemma}
\begin{proof}
The generalized central product of nilpotent groups is the
quotient of the direct product of finitely many nilpotent groups,
(which is again nilpotent, \cite{DS83} page 6) and another
nilpotent group. The quotient of nilpotent groups is nilpotent.
Therefore the generalized central product of finitely many
nilpotent groups is nilpotent.
\end{proof}
We note that the generalized central product of finitely many
solvable groups is solvable. For the case of abelian factors, the
generalized central product of an arbitrary number of abelian
groups is abelian.
\newline
Now we prove Theorem \ref{finiteNilpotentRS}:
\begin{proof} Suppose that $\{A_i| i \in I\}$ is an indexed family of
finitely generated nilpotent groups and let $G$ be the generalized
free product of the $A_i$, with amalgamating subgroup $C$:
\begin{eqnarray*}
G = \{{\prod_{i \in I}}^{\ast} A_i ; C \}.
\end{eqnarray*}
Let $S$ be the generalized central product of $A_i$ (see
Definition \ref{gen_cent_prod_def}):
\begin{eqnarray*}
S = {{\prod_{i \in I}}^{\times} A_i}/{gp(c \phi_i \; c^{-1}
\phi_j| i,j \in I,\; c \in C)}.
\end{eqnarray*}
$S$ contains a canonical homomorphism of $A_i$ for all $i \in I$.
Since these homomorphisms coincide on $C$, they can be extended to
a homomorphism $\mu$ from $G$ into $S$. Let $K$ be the kernel of
$\mu$. Since $S$ is nilpotent, 
it follows that $G/K$ is nilpotent. By Lemma
\ref{central_any_factor}, $\mu$ is one-to-one restricted to each
factor, i.e.
\begin{eqnarray*}
K \cap A_i =1 \text{ for all } i \in I.
\end{eqnarray*}
So, by the theorem of Hanna Neumann mentioned in Section
\ref{Hanna}, $K$ is free. Hence $G$ is a nilpotent extension of a
free group, so is also residually solvable.
\end{proof}
The following proposition is for the case of two finitely
generated solvable factors.
\begin{proposition} \label{central} The generalized free product of two finitely generated
solvable groups amalgamated by central subgroups, is a solvable
extension of a residually solvable group. Furthermore such groups
are residually solvable.
\end{proposition}
\begin{proof} Let $N$ be the normal closure of $B$ in $G$, i.e.
$N = {gp}_G (B)$ Since $A/C$ is solvable and $G/N \simeq A/C$ then
$G/N$ is a solvable group. Since $N$ is an iterated, infinite
untwisted double of copies of $B$, which therefore can be mapped
onto $B$ with free kernel, $N$ is residually solvable. So $G$ is a
solvable extension of a residually solvable group.
\end{proof}
For the case of abelian factors, we have the following corollary:
\begin{corollary} \label{infiniteabelianRS} The generalized free product an of
arbitrary number of abelian groups is an extension of a free group
by an abelian group. Further, such groups are residually solvable.
\end{corollary}
\section{Proof of Theorem \ref{abelianonefactor}}
Here we prove Theorem \ref{abelianonefactor}. Note that this is a
special instance of the case where the amalgamating subgroup is
central in only one of the factors.
\begin{proof} Let $A$ and $B$ be finitely generated groups, with $A$ torsion-free abelian and $B$ torsion-free nilpotent,
and let $G$ be their amalgamated product of them with amalgamating
subgroup $C$,
\begin{eqnarray*}
G= \{A \ast B;C_A = C_B\}.
\end{eqnarray*}
Then $C$ is a direct factor of a subgroup $A_1$ of $A$ of finite
index, i.e.
\begin{eqnarray*}
A_1 = C \times H \;\; \text{ where } [A : A_1] = n < \infty.
\end{eqnarray*}
There is an epimorphism from $G$ onto $A/{A_1} = \bigcup_{i =1}^n
a_i A_1$, where the $a_i$ are a distinct set of coset
representatives of $A_1$ in $A$. The kernel of this map is $K =
gp_G (B, A_1) = \{{\prod_{i = 1}^n}^{\ast} B^{a_i} \ast A_1; C\}$.
Now let $D$ be the normal closure of finitely many copies of $B$
in $K$:
\begin{eqnarray*}
D = gp_K (\bigcup_{i=1}^n B^{a_i}) = \{{{\prod^n_{i=1}}}^*
B^{a^i}; C\}.
\end{eqnarray*}
Since $[a_i, C]=1$, then $D$ is the generalized free product of a
finite number of doubles, and hence $D$ is residually solvable. So
$K$ is an extension of a residually solvable group by an abelian
group. Therefore $G$ is residually solvable-by-abelian-by-finite
abelian. Therefore $G$ is residually solvable.
\end{proof}
\subsection{An example where one of the factors is abelian}
Here we construct an example to illustrate the above theorem.
\begin{example} Let $A = gp(a,b,c)$ be a free abelian group, and $B
= gp(x, y, z)$ be a free nilpotent group of class $2$.
$C_A=gp(a^2, b)$ and $C_B = gp(y, z)$ are free abelian groups of
rank $2$. If we form
\begin{eqnarray*}
G=\{ A \ast B ; a^2 = x , b = z \},
\end{eqnarray*}
then $G$ is residually solvable. Note that $C_B$ is normal but not
central in $B$.
\end{example}
\begin{proof} The presentations of $A$, $B$, $C_A$ and $C_B$ are as
follows:
\begin{eqnarray*}
&& A = \langle a, b, c; [a,b], [b,c], [a,c] \rangle  \; \text{ free of rank $3$},\\
&& C_A = gp(a^2 , b) = \langle a^2, b; [a^2,b] \rangle  \; \text{ free abelian of rank $2$},\\
&& B = \langle x, y, z; [x,y]=z, [x,z], [y,z] \rangle  \; \text{ free nilpotent of class $2$ and of rank $3$},\\
&& C_B = gp(y ,z)=\langle y, z; [y, z] \rangle  \; \text{free
abelian of rank $2$}.
\end{eqnarray*}
Define an isomorphism $\phi$ which maps
\begin{eqnarray*}
a^2 \mapsto y \; \text{ and } b \mapsto z.
\end{eqnarray*}
Form the generalized free product of $A$ and $B$ identifying $C_A$
with $C_B$:
\begin{eqnarray*}
G = \{A \ast B; C_A = C_B\}.
\end{eqnarray*}
$C_B$ is normal in $B$ but is not central in $B$, and observe that $\xi B \cap C_B \not= \{1\}$.\\
However, $C_A$ is central and therefore normal in $A$. On the
other hand the quotient group $B/{C_B}$ is an infinite cyclic
group generated by $x C_B$,
\begin{eqnarray*}
B/{C_B} = gp(x C_B).
\end{eqnarray*}
Let $K$ be the normal closure of $A$ in $G$,
\begin{eqnarray*}
K= gp_G (A) = gp_B (A) = gp(A^{x^i}| i \in \mathbb{Z}).
\end{eqnarray*}
$G / K$ is an infinite cyclic group,
\begin{eqnarray*}
G / K = gp ( x K).
\end{eqnarray*}
Note that
\begin{eqnarray*}
a^{x^i} =b^i a^2 \;\; \text{ and } b^{x^i} =b \; \text{for all } i
\in \mathbb{Z}.
\end{eqnarray*}
Hence
\begin{eqnarray*}
A \cap A^{x^i} = C_A \;\;\; (\forall i \in \mathbb{Z})
\end{eqnarray*}
so defining
\begin{eqnarray*}
A_i := gp(A, A^{x^i})= \{ A \ast A^{x^i}; C_A \}
\end{eqnarray*}
we have
\begin{eqnarray*}
\forall w_i \in A_i \; \exists w_j \in A_j \text{ s.t. } w_i =
w_j^{x^{(i-j)}}
\end{eqnarray*}
Thus \begin{eqnarray*} K = \langle  A_i, x^i (i \in \mathbb{Z});
H_i = {H_j}^{x^{(i-j)}}(i,j \in \mathbb{Z}) \rangle.
\end{eqnarray*}
Then $K$ is an HNN-extension, where the base groups are the union
of $A_i \; (i \in {\mathbb Z})$, the stable letters are $x_i = x^i
(i \in \mathbb{Z})$, and the associated subgroups are $H_i$, where
$H_i$ is a subgroup of $A_i$ (for $i \in {\mathbb Z}$). In
particular by using the subgroup theorem, we can express $K$ as
the generalized free product of conjugates of free abelian groups.
\begin{eqnarray*}
K = \{{\prod_{i \in {\mathbb Z}}}^* A^{x^i}; C_A\}.
\end{eqnarray*}
Hence $K$ is an extension of a free group by an abelian group.
Consider the following short exact sequence:
\begin{eqnarray*}
1 \rightarrow K \rightarrow G \rightarrow B/{C_B} \simeq \mathbb{
Z } \rightarrow 1.
\end{eqnarray*}
$G/K$ is an infinite cyclic group, so therefore $G$ is an
extension of a residually solvable group by an infinite cyclic
group, and thus $G$ is residually solvable.
\end{proof}
\section{Proof of Theorem \ref{finind}}
We first give some background we need to prove Theorem
\ref{finind}.
\subsection{The Mal'cev completion and the fundamental theorem of
torsion-free nilpotent groups} Here we recall some definitions and
theorems that we use to prove Theorem \ref{finind}. A group is
called complete if, for every element $a$ of $G$ and natural
number $n$, the equation $x^n =a$ has at least one solution in
$G$, or, in other words, every root of every element of $G$
belongs to $G$ . There is a theorem that states that a finitely
generated nilpotent group is complete if and only if it contains
no proper subgroup of finite index; hence every group is contained
in a complete group. Let $\mathcal{D}$ denote the set of the class
of groups in which extraction of roots is always uniquely
possible. Let us term a minimal torsion-free nilpotent
$\mathcal{D}$-group, $m(G)$, containing a given torsion-free
nilpotent group $G$, a Mal'cev completion of $G$.

Let $G$ be a finitely generated torsion-free nilpotent group of
class $c$, let $g$ be an element of $G$ and let $n$ be a positive
integer. Then $G$ can be embedded in a finitely generated
torsion-free nilpotent group of class $c$ in which $g$ has an
$n$-th root. Therefore every finitely generated torsion-free
nilpotent group can be embedded in a nilpotent
$\mathcal{D}$-group.

We recall the Mal'cev Theorem: let $A$ be a torsion-free nilpotent
group, then $A$ can be embedded in a torsion-free nilpotent
$\mathcal{D}$-group. As a corollary we have: let $G$, $H$ be
torsion-free nilpotent groups, let $\phi$ be a homomorphism of $G$
into $H$ and let $m(G)$ and $m(H)$ be any Mal'cev completions of
$G$ and $H$; then $\phi$ can be extended uniquely to a
homomorphism $m(\phi)$ of $m(G)$ onto $m(H)$. Also if $m(G)$ and
$m'(G)$ are Mal'cev completions of $G$ then they are isomorphic.
We now state Mal'cev's fundamental theorem of torsion-free
nilpotent groups. If $G_1 ^*$ and $G_2 ^*$ are two completions of
$G$, then there exists an isomorphism between them that extends
the identity automorphism of $G$, and this isomorphism is unique
(for more details see \cite{GB71_}, \cite{AGK60}).
\subsection{Proof of Theorem \ref{finind}}
Note that here the condition of being torsion-free is necessary,
because of Mal'cev's Fundamental Theorem.
\begin{proof} Let $A$ and $B$ be two finitely generated torsion-free nilpotent
groups, and let $C$ be a proper subgroup of $A$ and $B$, such that
the index of $C$ in, say, $A$ is finite. We want to show that the
generalized free product of $A$ and $B$ with $C$ amalgamating,
\begin{eqnarray*}
G=\{ A \ast B; C \},
\end{eqnarray*}
is a finitely generated torsion-free nilpotent extension of a free
group. Let $m(A)$ and $m(B)$ be Mal'cev completions of $A$ and $B$
respectively. Since the index of $C$ in $A$ is finite, then $m(A)$
is also a Mal'cev completion of $C$. There exists a monomorphism
$\mu$ from $C$ into $B$. Then $\mu$ can be extended uniquely to a
homomorphism $m(\mu)$ of $m(A)$ onto $m(B)$. Now, there is a
homomorphism from $G$ into $m(B)$, which is consistent on $C$, and
injective when restricted to $A$ and $B$. Let $K$ be the kernel of
this homomorphism, then we have
\begin{eqnarray*}
K \cap A = \{1\} = K \cap B.
\end{eqnarray*}
Therefore $K$ is free by the result of Hanna Neumann mentioned in
Section \ref{Hanna}. Hence $G$ is a torsion-free nilpotent
extension of a free group.
\end{proof}
\section{Proof of Theorem \ref{finite_double}}
We first need to prove the following lemma which will be used in
proving Theorem \ref{finite_double}.
\begin{lemma} \label{double_key_lemma} If $A$ is a
group, $C$ is a subgroup of $A$, $\phi$ is an isomorphic mapping
of $A$ onto a group $B$, and $D$ is the amalgamated product of $A$
and $B$ amalgamating $C$ with $C \phi$, that is
\begin{eqnarray*}
D = \{ A \ast B ; C = C \phi \},
\end{eqnarray*}
then there is a homomorphism, $\psi$, from $D$ onto one of the
factors, with kernel:
\begin{eqnarray*}
\ker \psi = gp(a (a \phi)^{-1} | a \in A).
\end{eqnarray*}
Furthermore $\psi$ injects into each factor.
\end{lemma}
\begin{proof} Let $\alpha$ be the homomorphism from $A$ onto itself, and $\beta$ be
the homomorphism from $B$ onto the inverse of the isomorphic copy
of $A$, i.e. $\beta = {\phi}^{-1}$. These homomorphisms can be
extended to a homomorphism from $D$ onto $A$, (\cite{GB93} page
103, \cite{BN49}). Since $\psi(C a b) = C a (b \phi^{-1})$ it
follows easily $K = \ker \psi = gp (a (a \phi^{-1})| a \in A).$ By
the way that $\alpha$ and $\beta$ are defined, it follows that
this homomorphism is one-to-one restricted to either $A$ or $B$.
\end{proof}
Now we are ready to prove Theorem \ref{finite_double}.
\begin{proof} For each $i \in I$ let $\phi : G \rightarrow A_i$ be an epimorphism, and let $K$ be the kernel of $\phi$.
$K$ is free, since $\phi$ restricted to each factor is injective,
and
\begin{eqnarray*}
A_i \cap K =\{1\} \;\;(\forall i \in I).
\end{eqnarray*}
Therefore $G$ is free-by-nilpotent.
\end{proof}
\section{Proof of Proposition \ref{nil_neg}}
\begin{proof}
Let $A = gp(a,x)$ be a free nilpotent group of class $2$, and $B =
gp(b,y)$ be a free nilpotent group of class $3$. Then $C_A=gp(a,
a^x)$ and $C_B = gp( b, [b, b ^y])$ are free abelian groups of
rank $2$. If we form
\begin{eqnarray*}
  G & = & \{ A \ast B ; C_A = C_B \}\\
    & = & \{ A \ast B ; a=b, a^x = [b , b^y] \},
\end{eqnarray*}
then
\begin{eqnarray*}
  a & = & [a, a^y]^x = [a^x , a^{y^x}],
\end{eqnarray*}
so $a$ lies in every term of the derived series of $G$. Hence $G$
is not residually solvable. Now put $N={gp}_{G} (a)$ which is
perfect. Since a poly-residually solvable groups can not have a
perfect subgroup, $G$ is not poly-residually solvable.
\end{proof}
Note that this contrasts with the case where residual finiteness
of the groups has been studied, where Baumslag proved that
generalized free products of finitely generated torsion-free
nilpotent groups in general are meta-residually solvable.

We have shown that if the amalgamating subgroup of the generalized
free product of torsion-free nilpotent groups is cyclic
(Proposition \ref{Cyclic_amalgam}) or central in both factors
(Theorem \ref{central}) then these groups are residually solvable.
One may conjecture that it would be the case when the amalgamating
subgroup is abelian. However this is false by the example of
Baumslag as mentioned above.

This can be easily shown using Baumslag's example of an
amalgamated product, which can be modified to give an example
where the factors are finite, but the product is not residually
solvable.

\section{Proof of Theorem \ref{Poly-RS}}
\label{Poly_RS} Let us recall the definition of a poly-property.
Let $X$ be a property (or class) of groups. A finite normal series
\begin{eqnarray*}
1 = G_0 \leqslant G_{1} \leqslant \cdots \leqslant G_l = G
\;\;\;\;\; (1)
\end{eqnarray*}
of $G$ is termed a poly-$X$ series for $G$ if
\begin{eqnarray*}
{G_{i+1}}/{G_i} \in X \text{ for } i=0, \cdots, l-1.
\end{eqnarray*}
A group $G$ is termed poly-$X$ if it has at least one poly-$X$
series. The length of the series (1) is $l$.

The proof of Theorem \ref{Poly-RS} uses an upper central
filtration approach:
\begin{proof}
Let $G=\{A \ast B; C_A=C_B\}$ be the generalized free product of
finitely generated nilpotent groups $A$ and $B$. Assume that
$\xi_i A$ and $\xi_i B$ be $i$-th terms of the upper central
series of $A$ and $B$ respectively and $(\xi_i A), (\xi_i B)$ (for
$i$ is bounded by the maximum class of nilpotency of $A$ and $B$)
be solvable filtrations of $A$ and $B$ such that they are also
compatible filtrations. We want to show that $G$ has an invariant
series
\begin{eqnarray*}
1=G_0 \leq G_1 \leq \cdots \leq G_k = G \;\; (k < \infty)
\end{eqnarray*}
such that ${G_{i+1}}/{G_i}$ is residually solvable (for $i=0,
\cdots, k-1$). Map
\begin{eqnarray*}
G \rightarrow G/{gp_G (C_A \cap \xi_1 A)}.
\end{eqnarray*}
Note that, since $\xi G = \xi A \cap \xi B \cap C_A$, then
\begin{eqnarray*}
gp_G (C_A \cap \xi_1 A) < \xi_1 G.
\end{eqnarray*}
\begin{eqnarray*}
G/{gp_G (C_A \cap \xi_1 A)} \simeq &&\{A/{gp_G (C_A \cap \xi_1 A)} \ast B/{gp_G (C_A \cap \xi_1 A)};\\
                                   && {C_A}/{gp_G (C_A \cap \xi_1 A)} = {C_B}/{gp_G (C_A \cap \xi_1 A)}\}\\
G/{gp_G (C_A \cap \xi_1 A)} \simeq &&\{ A/{C_A \cap \xi_1 A} \ast B/{C_B \cap \xi_1 B};\\
                                   && {C_A}/{C_A \cap \xi_1 A} = {C_B}/{C_B \cap \xi_1 B)}\}.
\end{eqnarray*}
One can check that this is residually solvable.
\begin{eqnarray*}
&& {G/{gp_G (C_A \cap \xi_1 A)}}/{gp_G ({\xi_1 A}/{C_A \cap \xi_1
A},
{\xi_1 B}/{C_B \cap \xi B})} \simeq \\
&& \{ A/{\xi_1 A} \ast B/{\xi_1 B}; {C_A \xi_1 A}/{\xi_1 A} = {C_B
\xi_1 B}/{\xi_1 B}\}.
\end{eqnarray*}
Inductively we can show that for $i$ each quotient is residually
solvable and this completes the proof of the theorem.
\end{proof}
\bibliographystyle{amsplain}
\bibliography{XBib}
\end{document}